\documentclass[article,12pt]{article}
\usepackage{amsmath, latexsym, amsfonts, amssymb, amsthm, amscd}

\newtheorem{theorem}{Theorem}
\newtheorem{coro}{Corollary}
\newtheorem{prop}{Proposition}

\newtheorem{lemma}{Lemma}

\newcommand{\ind}{\mbox{\rm 1\hspace{-0.04in}I}}

\newcommand{\ed}{\stackrel{(d)}{=}}

\DeclareMathOperator{\e}{\mathbb{E}}
\DeclareMathOperator{\p}{\mathbb{P}}

\DeclareMathOperator{\R}{\mathbb{R}}

\numberwithin{equation}{section}

\begin{document}


\title{Conditioned stable L\'evy processes and\\ Lamperti
representation.}

\maketitle

\begin{center}
{\large M.E. Caballero\footnote{Instituto de Matem\'aticas,
Universidad Nacional Aut\'onoma de M\'exico, {\sc Mexico 04510
DF.} E-mail: emilia@servidor.unam.mx} and L.
Chaumont\footnote{Laboratoire de Probabilit\'es et Mod\`eles
Al\'eatoires, Universit\'e Pierre et Marie Curie, 4, Place Jussieu
- 75252 {\sc Paris Cedex 05.} E-mail: chaumont@ccr.jussieu.fr}}
\end{center}
\vspace{0.2in}

\begin{abstract}
By killing a stable L\'evy process when it leaves the positive
half-line, or by conditioning it to stay positive, or by
conditioning it to hit 0 continuously, we obtain three different
positive self-similar Markov processes which illustrate the three
classes described by Lamperti \cite{La}. For each of these
processes, we compute explicitly the infinitesimal generator from
which we deduce the characteristics of the underlying L\'evy process
in the Lamperti representation. The proof of this result bears on
the behaviour at time 0 of stable L\'evy processes before their
first passage time across level 0 which we describe here. As an
application, we give the law of the minimum before an independent
exponential time of a certain class of L\'evy processes. It provides
the explicit form of the spacial Wiener-Hopf factor at a particular
point and the value of the ruin probability for this class of
L\'evy processes.\\

\noindent {\sc Key words and phrases}: Positive self-similar Markov
processes, Lamperti representation, infinitesimal generator, stable
L\'evy processes conditioning to stay positive,
stable L\'evy processes conditioning to hit 0 continuously.\\

\noindent MSC 2000 subject classifications: 60 G 18, 60 G 51, 60 B
52.
\end{abstract}

\vspace{0.5in}

\section{Introduction and preliminary results}

The stochastic processes which are considered in this work take
their values in the Skohorod's space $\mathcal{D}$ of c\`adl\`ag
trajectories. We define this set as follows: $\Delta:=+\infty$ being
the cemetery point, a function
$\omega:[0,\infty)\rightarrow\R\cup\Delta$ belongs to $\mathcal{D}$
if and only if:
\begin{itemize}
\item[--] For all $t\ge\zeta(\omega)$, $\omega_t=\Delta$, where
$\zeta(\omega):=\inf\{t:\omega_t=\Delta\}$ is the lifetime of
$\omega\in\mathcal{D}$ and $\inf\emptyset=+\infty$.
\item[--]
For all $t\ge0$, $\lim_{s\downarrow t}\omega_s=\omega_t$ and for all
$t\in(0,\zeta(\omega))$, $\lim_{s\uparrow t}\omega_s:=w_{t-}$ is a
finite real value.
\end{itemize} The space $\mathcal{D}$ is endowed with the Skohorod's
$J_1$ topology. We denote by $X:\mathcal{D}\rightarrow\mathcal{D}$
the canonical process of the coordinates and by $({\cal F}_t)$ the
natural Borel filtration generated by $X$, i.e. ${\cal
F}_t=\sigma(X_s,\,s\le t)$. A probability measure $P_x$  on
$\mathcal{D}$ is the law of a L\'evy process if $(X,P_x)$ starts
from $x$, i.e.~$P_x(X_0=x)=1$ and has independent and homogeneous
increments.  Note that $(X,P_x)=(x+X,P_0)$ and that the lifetime of
$(X,P_x)$ is either a.s.~infinite or  a.s.~finite. It is well know
that for any L\'evy process $(X,P_x)$ with finite lifetime
$\zeta(X)$, there exists a L\'evy process $(X',P_x)$ with infinite
lifetime, such that under $P_x$ the random variable $\zeta(X)$ is
exponentially distributed and independent of $X'$ and $P_x$-a.s.,
$X_t=X'_t$, if $t<\zeta(X)$. Furthermore the
parameter of the law of of $\zeta(X)$ under $P_x$ does not depend on $x$.\\

\noindent An $\R_+$-valued self-similar Markov process $(X,\p_x)$,
$x>0$ is a strong Markov process with values in the space
$\mathcal{D}$, which fulfills a scaling property, i.e.~there
exists a constant $\alpha > 0$ such that:
\begin{equation}\label{scale}
\mbox{\it The law of $\;(kX_{k^{-\alpha}t},\,t\ge0)$ under $\p_x$ is
$\p_{kx}$.}
\end{equation}
We will call these processes pssMp for short.  They are much
involved in many areas of probability theory. For instance, the
continuous state branching process obtained as the weak limit of a
re-scaled discrete branching process is a pssMp which is associated
to a self-similar L\'evy tree, see \cite{DL}. These processes also
appear in fragmentation theory~; the mass process of a self-similar
fragmentation process is itself a pssMp, \cite{Be2}. The pssMp that
we are going to study here have recently been obtained in \cite{CC}
as limits of re-scaled random walks whose law are in the domain of
attraction of a stable law, after they are conditioned to stay
positive or conditioned to hit 0 at a finite time, see sections
\ref{pos} and \ref{zero} below.

According to Lamperti \cite{La}, the set of pssMp splits into three
exhaustive classes which can be distinguished from each other by
comparing their values at their first hitting time of 0, i.e.:
\[S = \inf \{t>0 : X_t =0\}. \]
This classification may be summerized as follows:
\begin{itemize}
\item $\mathcal{C}_1$ is the class of pssMp such that $S =
+\infty $, $\p_x$-a.s. for all starting points $x>0$. \item
$\mathcal{C}_2$ is the class of those for which $S< +\infty$ and
$X_{S-} = 0$, $\p_x$-a.s. for all starting points $x>0$. Processes
of this class hit the level 0 in a continuous way. \item
$\mathcal{C}_3$ is that of those for which $S < +\infty$ and $X_{S-}
> 0$, $\p_x$-a.s. for all starting points $x>0$. In that case, the
process hits 0 by a negative jump.
\end{itemize}
The main result of \cite{La} asserts that any pssMp up to its first
hitting time of 0 may be expressed as the exponential of a L\'evy
process, time changed by the inverse of its exponential functional.
Then the underlying L\'evy process in the so-called Lamperti
representation of $(X,\p_x)$ fulfils specific features depending on
the class to which $(X,\p_x)$ belongs. More formally, let $(X,\p_x)$
be a pssMp starting from $x>0$, and write the canonical process $X$
in the following form:
\begin{equation}\label{lamp}
\mbox{$X_t=x\exp\xi_{\tau(tx^{-\alpha})},\;\;\;0\le t<S\,, $}
\end{equation}
where for $t<S$,  $\tau(t) = \inf \{s\geq 0 : \int_0^s
\exp\alpha\xi_u du \geq t\}$. Then under $\p_x$, $\xi=(\xi_t,\;t\geq
0)$ is a L\'evy process started from $0$ which law does not depend
on $x>0$ and such that
\begin{itemize}
\item if  $(X,\p_x)\in\mathcal{C}_1$,  then $\zeta(\xi)=+\infty$
and $\limsup_{t\rightarrow+\infty}\xi_t=+\infty$, $\p_x$-a.s.

\item  if $(X,\p_x)\in\mathcal{C}_2$, then $\zeta(\xi)=+\infty$
and $\lim_{t\to\infty} \xi_t = -\infty$, $\p_x$-a.s.

\item if $(X,\p_x)\in\mathcal{C}_2$, then $\zeta(\xi)<+\infty$,
$\p_x$-a.s.
\end{itemize}
Note that for any $t < \int_0^\infty\exp(\alpha\xi_s)\,ds$,
$$\tau(t)=\int_0^{x^\alpha t}\frac{ds}{X_s^\alpha}\,,\;\;\;
\p_x-\mbox{a.s.}$$ so that (\ref{lamp}) is revertible and yields a
one to one relation between the class of pssMp's killed at time $S$
and the one of L\'evy processes.

Now we recall another important result of Lamperti \cite{La} which
gives the explicit form of the generator of any pssMp in terms of
this of its underlying L\'evy process. Let $(X,\p_x)$ and $\xi$ be
any such processes related as in (\ref{lamp}). We will denote by
$\mathcal{K}$ and $\mathcal{L}$ their respective generators and by
$\mathfrak{D}_{\mathcal{K}}$ and $\mathfrak{D}_{\mathcal{L}}$ the
domains of $\mathcal{K}$ and $\mathcal{L}$. Then recall that
$\mathfrak{D}_{\mathcal{L}}$ contains all the functions with
continuous second derivatives in $[-\infty, +\infty]$ and  if
$\tilde{f}$ is such a function then $\mathcal{L}$ is of the form:
\begin{equation}\label{lk}
\mathcal{L} \tilde{f}(x) = a \tilde{f}^{\prime}(x) +
\frac\sigma2\tilde{f}''(x)+\int_{\R} [\tilde{f}(x+y) - \tilde{f}(x)
- \tilde{f}^{\prime}(x)l(y)]\Pi(dy) - k
\tilde{f}(x)\,,\end{equation} for $x\in\R$, where $a\in\R$,
$\sigma>0$. The measure $\Pi(dx)$ is the L\'evy measure of $\xi$ on
$\R$; it verifies $\Pi(\{0\})=0$ and
$\int(1\wedge|x|^2)\,\Pi(dx)<\infty$. The function $l(\cdot)$ is a
bounded Borel function such that $l(y)\sim y$ as $y\rightarrow0$.
The last term $k\ge0$ corresponds to the killing rate of $\xi$, that
is the parameter of $\zeta(\xi)$, ($k=0$ if $\xi$ has infinite
lifetime). {\it It is important to note that in the expression
$(\ref{lk})$, the choice of the function $l(\cdot)$ is arbitrary and
the coefficient $a$ is the only one which depends on this choice.}

Theorem 6.1 of Lamperti \cite{La} may be stated as follows:
\begin{theorem}[Lamperti \cite{La}]\label{thlamp} If $f :[0,+\infty] \to \R$
is such that $f,xf^{\prime}, x^2f^{\prime\prime}$ are continuous in
$[0,+\infty]$, then they belong to the domain
$\mathfrak{D}_{\mathcal{K}}$ of the infinitesimal generator of
$(X,\p_x)$ which has the form
\begin{eqnarray*} \mathcal{K}f(x) &=&
\frac{1}{x^{\alpha}}\int_{\R^+}[f(ux)-f(x)-f^{\prime}(x)l(\log u)
]\Theta(du)\\&&\qquad\qquad\qquad a x^{1-\alpha}f^{\prime}(x)
+\frac{\sigma}2 x^{2-\alpha}f''(x)- k x^{-\alpha}f(x),
\end{eqnarray*}
for $x>0$, where $\Theta(du) = \Pi (d u)\circ\log u$, for $u>0$.
This expression determines the law of the process $(X_t,\,0\le t\le
S)$ under $\p_x$.
\end{theorem}

To present the results of this paper, let us first consider two
examples in the continuous case. The first one is when $(X,\p_x)$ is
the standard real Brownian motion absorbed at level 0. The process
$(X,\p_x)$ is a pssMp which belongs to the class ${\cal C}_2$, with
index $\alpha=2$ and it is well known (see for instance \cite{CPY})
that its associated L\'evy process in the Lamperti representation
(\ref{lamp}) is given by $\xi=(B_t-t/2,\,t\ge0)$, where $B$ is a
standard Brownian motion. The second example is when $(X,\p_x)$ is
the Brownian motion conditioned to stay positive. This process
corresponds to the three dimensional Bessel process, i.e.~the norm
of a three dimensional Brownian motion. Then, $(X,\p_x)$ is a pssMp
which belongs to the class ${\cal C}_1$,  with index $\alpha=2$ and
the underlying L\'evy process is given by $\xi=(B_t+t/2,\,t\ge0)$.

Similarly, it is possible to obtain pssMp's from any stable L\'evy
process $(X,P_x)$ with index $\alpha\in(0,2)$, throughout the same operations.
More precisely, by killing $(X,P_x)$ when it enters into the
negative halfline, i.e.
\[X_t\ind_{\{t<T\}}\,,\;\;\mbox{with}\;\;T=\inf\{t\ge0:X_t\le0\}\,,\]
 one obtains a pssMp $(X,\p_x)$ which belongs to the class
${\cal C}_2$ or ${\cal C}_3$ according as $(X,P_x)$ has negative jumps or
not. Also by conditioning a stable L\'evy process to stay positive, i.e.
\[\p_x=\lim_{t\rightarrow+\infty}P_x(\,\cdot\,|\,T>t)\,,\;\;x>0\,,\]
one obtains a pssMp $(X,\p_x)$ belonging to ${\cal C}_1$. One may
also give a sense to the conditioning  to hit 0 continuously; such
processes belong to ${\cal C}_3$. The main goal of this paper, is to
identify the underlying L\'evy process in the Lamperti
representation for each of these processes by computing their
infinitesimal generators and using Lamperti's result recalled above.
This will be done in section 3. In section 4, we deduce from the
results of section 3.1, the law of the minimum before an independent
exponential time for an important class of L\'evy processes. It
gives the expression of the Wiener-Hopf factor of these L\'evy
processes at a particular point, i.e. the law of $\inf_{s\le{\bf
e}(k)}{\xi_s}$, where $\xi$ is a L\'evy process which
characteristics are described in Corollary \ref{cor1} and ${\bf
e}(k)$ is an independent random variable with a special parameter
$k$. We also find the law of the overall minimum for another class
of L\'evy processes whose law is given by Corollary 2. This
calculation is equivalent to the problem of finding the explicit
form of the corresponding ruin probability which has recently been
studied for other classes of L\'evy processes by Lewis and Mordecki
\cite{LM}. The next section is devoted to further preliminary
results, the main of which having some interest in its own,
independently of the rest of the paper. It extends a result of
Bingham \cite{Bi} and Rivero \cite{Ri} which describes the
asymptotic behaviour as $t$ goes to 0 of $P_x(T\le t)$, that is the
small tail of first passage times of stable L\'evy processes.

\section{Small tail of first passage times of stable L\'evy
processes}\label{prelim}

In all the sequel of this paper, $(X,P_x)$ will be a stable L\'evy
process with index $\alpha\in(0,2)$, starting at $x\in\R$. Since stable
L\'evy processes have infinite lifetime,
the characteristic exponent of $(X,P_x)$ is defined by $E_0[\exp(i\lambda
X_t)]=\exp[t\psi(\lambda)]$, $t\ge0$, $\lambda\in\R$, where
\begin{equation}\label{2938}
\psi(\lambda)=ia\lambda+
\int_{\R}(e^{i\lambda y}-1-i\lambda y\ind_{\{|y|<1\}})\,\nu(y)\,dy\,.
\end{equation}
The density of the L\'evy measure is
\begin{equation}\label{1894}
\nu(y) = c_+ y^{-\alpha-1}\textbf{1}_{\{y>0\}}
        + c_- \vert y \vert ^{-\alpha-1}\textbf{1}_{\{y<0\}}\,,\end{equation}
where $c_+$ and $c_-$ are two nonnegative constants such that
$c_++c_->0$. Note also that the constant $a$ is related to $c_+$,
$c_-$ and $\alpha$ as follows: $a = \frac{c_+-c_-}{1-\alpha}$,
$\alpha\neq 1$. In the case where $\alpha=1$, the process $(X,P_x)$
will be supposed to be a symmetric Cauchy process, so we have
$c_+=c_-$ and $a=0$. We suppose
moreover that neither $(X,P_x)$ nor $(-X,P_x)$ is a subordinator.\\

The main result of this section concerns the asymptotic behaviour as
$t\downarrow0$ of
\[E_x(f(X_t)\ind_{\{T\le t,\,X_t\in(0,\infty)\}})\,,\;\;\mbox{with}\;\;
T=\inf\{t:X_t\le0\}\,,\] where $f$ is a bounded and continuous
function.  This result will be used to compute the infinitesimal
generator of the killed stable L\'evy processes. We denote by
$P_x(\,\cdot\,|\,X_t=y)$ a regular version of the law of the bridge
of the L\'evy process $(X,P_x)$ from $x$ to $y$, with length $t$.
Let $p_s(z)$, $s\ge0$, $z\in\R$ be the density of the semigroup of
$(X,P_x)$, then for all $s\in[0,t)$, this law is defined on
$\mathcal{F}_s$ by
\begin{equation}\label{357}P_x(A\,|\,X_t=y)=E_x\left(\ind_{A}
\frac{p_{t-s}(y-X_{s})}
{p_t(y-x)}\right)\,,\;\;\;A\in\mathcal{F}_s\,.\end{equation} See
\cite{FPY} for a complete account on bridges of Markov processes.

Now let us recall some classical properties of densities of stable
laws which may be found in \cite{Zo} or \cite{Sa}, Chap.~3.14. When
the corresponding L\'evy measure is not concentrated on either
$(-\infty,0]$ or $[0,\infty)$, there are constants $C_1,C_2>0$, such
that
\begin{equation}\label{238}
p_1(z)\sim C_1|z|^{-\alpha-1}\,,\;\;\;\mbox{as
$z\rightarrow-\infty\;\;$ and}\;\;\;\;p_1(z)\sim
C_2z^{-\alpha-1}\,,\;\;\;\mbox{as $z\rightarrow+\infty$.}
\end{equation}
If the L\'evy measure is concentrated on $(-\infty,0]$, then there
are constants $C_3,C_4>0$ such that
\begin{equation}\label{239}
p_1(z)\sim C_3|z|^{-\alpha-1}\,,\;\;\;\mbox{as
$z\rightarrow-\infty\;\;$ and}\;\;\;\;p_1(z)\sim C_4
x^{2-\alpha}e^{-x}\,,\;\;\;\mbox{as $z\rightarrow+\infty$},
\end{equation}
where $x=(\alpha-1)(z/\alpha)^{\alpha/(\alpha-1)}$. Note that in
this second case, we have necessarily $1<\alpha<2$ since we have
implicitly excluded subordinators of our study.

Our first lemma expresses the intuitive fact that the amplitude of
a bridge from $x$ to $y$ of $(X,P_x)$ tends to $|y-x|$ as its
length goes to 0. It might be established more directly from a
suitable estimation of the join law of $(X_t,\underline{X}_t)$,
under $P_x$, where $\underline{X}_t:=\inf_{s\le t}X_s$, however we
have not found any such result in the literature.

\begin{lemma}\label{lem1} For all $x,y>0$,
\[\lim_{t\rightarrow0}P_x(T\le t\,|\,X_t=y)=0\,.\]
\end{lemma}
\begin{proof} First let $t>0$ and decompose the term
of the statement as
\begin{equation}\label{439}P_x(T\le t\,|\,X_t=y)=P_x(T\le t/2\,|\,X_t=y)+
P_x(T\in(t/2,t]\,|\,X_t=y)\,.\end{equation}
To prove the result, it is enough to show that the first term in (\ref{439})
converges to 0 as $t$ tends to 0. Indeed, let $(X,\hat{P}_x):=(-X,P_x)$ be
the dual L\'evy process, then the
following identity in law between the bridge and its time reversed
version is well known, see \cite{FPY} for instance:
\begin{equation}\label{return}
\big((X_{(t-s)-},\,0\le s\le
t),P_x(\,\cdot\,|\,X_t=y)\big)=\big((X_{s},\,0\le s\le
t),\hat{P}_y(\,\cdot\,|\,X_t=x)\big)\,.
\end{equation}
(We have set $X_{0-}=0$.) Then we observe the inequality:
\[P_x(T\in(t/2,t]\,|\,X_t=y)\le\hat{P}_y(T\le t/2\,|\,X_t=x)\,.\]
If the first term of (\ref{439}) converges to 0 in any case, then by
applying the result to the bridge of the dual process and the above
inequality, we show that the second term of (\ref{439}) converges
also to 0.

Now, let us prove that the first term of (\ref{439}) converges to 0 as $t$ goes
to $0$. Recall that
$\underline{X}_t:=\inf_{s\le t}X_s$. From (\ref{357}) the
first term is
\begin{eqnarray}P_x(T\le t/2\,|\,X_t=y)&=&E_x\left(\ind_{\{T\le
t/2\}}\frac{p_{t/2}(y-X_{t/2})}{p_t(y-x)}\right)\nonumber\\
&=&E_0\left(\ind_{\{\underline{X}_{t/2}\le-x\}}
\frac{p_{1}(2^{1/\alpha}t^{-1/\alpha}[y-x-X_{t/2}])}
{2^{-1/\alpha}p_1(t^{-1/\alpha}[y-x])}\right)\,,\label{901}
\end{eqnarray}
where the second identity follows from the fact that
$p_t(z)=t^{-1/\alpha}p_1(t^{-1/\alpha}z)$, for all $t>0$.

From classical properties of stable L\'evy processes, we have
$P_0(\underline{X}_{t/2}\le-x)\rightarrow0$ as $t\rightarrow0$ and
$p_1(0)>0$. Therefore, if $x=y$, then since $z\mapsto p_1(z)$ is
bounded on $\R$, we see that the right hand side of (\ref{901})
tends to 0 as $t$ goes to $0$. So the lemma is proved when $x=y$.

Set $q=1-2^{-1/(2\alpha)}$ and suppose that $y> x$, then again we
develop the right hand side of (\ref{901}) as the sum:
\begin{eqnarray}&&
E_{0}\left(\ind_{\{\underline{X}_{t/2}\le-x,\,X_{t/2}\le q(y-x)\}}
\frac{p_{1}(2^{1/\alpha}t^{-1/\alpha}[y-x-X_{t/2}])}
{2^{-1/\alpha}p_1(t^{-1/\alpha}[y-x])}
\right)\nonumber\\&&\;\;\;\;\;\;\;\;\;\;\;\;+E_{0}
\left(\ind_{\{\underline{X}_{t/2}\le -x,\,X_{t/2}\ge q(y-x)\}}
\frac{p_{1}(2^{1/\alpha}t^{-1/\alpha}[y-x-X_{t/2}])}
{2^{-1/\alpha}p_1(t^{-1/\alpha}[y-x])}\right).\label{862}
\end{eqnarray}
Note that on the event $\{X_{t/2}\le q(y-x)\}$, we have
$$2^{1/\alpha}t^{-1/\alpha}[y-x-X_{t/2}]\ge t^{-1/\alpha}[y-x]>0.$$ So,
from (\ref{238}) and (\ref{239}), there is a time $t_1$ and a finite
constant $c_1$ (both non random) such that for all $0<t\le t_1$, on
the event $\{X_{t/2}\le q(y-x)\}$ we have
\begin{equation}
\frac{p_{1}(2^{1/\alpha}t^{-1/\alpha}[y-x-X_{t/2}])}
{p_1(t^{-1/\alpha}[y-x])}\le c_1\,.\end{equation} Hence from
Lebesgue theorem of dominated convergence the first term in
(\ref{862}) tends to 0 at $t$ goes to 0. Now call
$\hat{p}_t(z):=p_t(-z)$ the semigroup of the dual process
$(X,\hat{P}_0)$. Since bridges of L\'evy processes have no fixed
discontinuities, see \cite{FPY}, from (\ref{return}), the second
term in (\ref{862}) may be written as
\begin{eqnarray*}
&&E_{0}\left(\ind_{\{\underline{X}_{t/2}\le-x,\,X_{t/2}\ge q(y-x)\}}
\frac{p_{1}(2^{1/\alpha}t^{-1/\alpha}[y-x-X_{t/2}])}
{2^{-1/\alpha}p_1(t^{-1/\alpha}[y-x])}\right)\\
&=&P_x(\underline{X}_{t/2}\le0,\,X_{t/2}\ge q(y-x)+x\,|\,X_t=y)\\&=&
\hat{P}_y(\underline{X}_{t/2}\le 0,\,X_{t/2}\ge
q(y-x)+x\,|\,X_t=x)\\&=&
\hat{E}_{0}\left(\ind_{\{\underline{X}_{t/2}\le 0,\,X_{t/2}\ge
(1-q)(x-y)\}}
\frac{\hat{p}_{1}(2^{1/\alpha}t^{-1/\alpha}[x-y-X_{t/2}])}
{2^{-1/\alpha}\hat{p}_1(t^{-1/\alpha}[x-y])}\right)\,.
\end{eqnarray*}
On the event $\{X_{t/2}\ge
(1-q)(x-y)\}$, we have
$$2^{1/\alpha}t^{-1/\alpha}[x-y-X_{t/2}]\le (2^{1/\alpha}-2^{1/(2\alpha)})
t^{-1/\alpha}[x-y]<0.$$ If $(X,P_x)$ has positive jumps, then from
(\ref{239}), $\hat{p}_1(z)\sim C_3 |z|^{-\alpha-1}$ as
$z\rightarrow-\infty$, thus there is a time $t_2$ and a finite
constant $c_2$ (both non random) such that for all $0<t\le t_2$, on
the event $\{X_{t/2}\ge (1-q)(x-y)\}$ we have
\begin{equation}
\frac{\hat{p}_{1}(2^{1/\alpha}t^{-1/\alpha}[x-y-X_{t/2}])}
{\hat{p}_1(t^{-1/\alpha}[x-y])}\le c_2\,,\end{equation} hence, again
the second term in (\ref{862}) tends to 0 at $t$ goes to 0.

So we have proved the lemma when $y>x$ and $(X,P_x)$ has positive
jumps. By a time reversal argument, it is easy to see that the same
result holds when $y<x$ and $(X,P_x)$ has negative jumps. It remains
to show the result when $y<x$ and $(X,P_x)$ has no negative jumps.
In this case, put $T_y=\inf\{t:X_t=y\}$, then from the Markov
property applied at time $T_y$, we have
\begin{equation}\label{100}P_x(T\le t\,|\,X_t=y)=
\int_0^t P_x(T_y\in ds\,|\,X_t=y)P_y(T\le
t-s\,|\,X_{t-s}=y)\,.\end{equation} But we already proved above that
$P_y(T\le t-s\,|\,X_{t-s}=y)$ tends to 0 as $t-s$ goes to 0. This
ends the proof of the lemma.
\end{proof}

\noindent Recall that the characteristic exponent of $(X,\p)$ may also be written
in the following form for $\alpha\in(0,1)\cup(1,2)$,
\begin{equation}\label{2238}
E_0[\exp(i\lambda X_t)]=\exp[-ct|\lambda|^\alpha(1-i\beta
\mbox{sgn}(\lambda)\tan(\pi\alpha/2))]\,,\;\;\; \lambda\in\R\,,\end{equation}
where
\[c=(c_++c_-)\Gamma(-\alpha)\cos\frac{\pi\alpha}2\,\;\;\mbox{and}\,
\;\;\beta=(c_+-c_-)/(c_++c_-)\,,\]
see for instance Sato \cite{Sa}, Theorem 14.10 and its proof p.83--85.
It has been proved by Bingham \cite{Bi}, Proposition 3.b and Theorem 4.b, and
Rivero \cite{Ri}, section 2.3 that
\begin{equation}\label{437}
\lim_{t\downarrow0}\frac1tP_x(T\le t)=
\frac{k}{x^{\alpha}}\,,\end{equation} where the constant $k$ is
explicitly computed in \cite{Bi} and is given by:
\begin{equation}\label{1257} k=
c(1+\beta^2\tan^2(\pi\alpha/2))^{1/2}\Gamma(\alpha)
\sin(\pi\alpha\rho)/\pi\,.\end{equation}
By definition, $\rho:=P_0(X_1<0)$ and it is well known that this rate has
the expression
\[\rho=\frac12-(\pi\alpha)^{-1}\arctan(\beta\tan(\pi\alpha/2))\,.\]
Note that we always have $\alpha\rho\le1$. Moreover, we easily check
that $(X,P_x)$ has no negative jumps if and only if one of the three
following conditions holds
\[c_-=0\Leftrightarrow\beta=1\Leftrightarrow\alpha\rho=1\,.\]
For $\alpha=1$, the expressions (\ref{2238}) and (\ref{1257}) are
reduced to $E_0[\exp(i\lambda X_t)]=\exp(-c_+\pi t|\lambda|)$ and
$k=c_+\,(=c_-)$, respectively (although the value of $k$ in this
case is ambiguous in \cite{Bi}, it will be confirmed in section
\ref{pos}).

Then in section \ref{pos} we will provide another means to compute
the expression of the constant $k$, see formula (\ref{constk}). Note
also that Rivero's result \cite{Ri} concerns the more general
setting of positive self-similar Markov processes. Besides, in the
case where $(X,P_x)$ has no negative jumps, we have $k=0$ but
Proposition 3.b of \cite{Bi} gives an explicit form of the
asymptotic behaviour of $P_x(T<t)$, as $t\downarrow0$. The next
theorem completes Bingham and Rivero's result.

\begin{theorem}\label{lem2} For all $x>0$, and all bounded,
continuous function $f:\R\rightarrow\R$,
\[\lim_{t\rightarrow0}\frac1tE_x(f(X_t)\ind_{\{T\le t,\,X_t\in(0,\infty)\}})=
\frac{f(x)}{x^\alpha}\left(k-\frac{c_-}{\alpha}\right)\,,\]
where $c_-$ and $k$ are respectively defined in $(\ref{1894})$ and $(\ref{1257})$.
\end{theorem}
\begin{proof}
Let $\delta\in(0,x)$ and  define
$I_{\delta,x}:=[x-\delta,x+\delta]$. Let also $f$ be a bounded and
continuous function and write: \begin{eqnarray*} \frac1t
E_x(f(X_t)\ind_{\{T\le t,\,X_t\in(0,\infty)\}})&=&\frac1t
E_x(f(X_t)\ind_{\{T\le t,\,X_t\in
I_{\delta,x}\cap(0,\infty)\}})\\&&\qquad+\frac1t
E_x(f(X_t)\ind_{\{T\le t,\,X_t\in
I_{\delta,x}^c\cap(0,\infty)\}})\,.\end{eqnarray*} Then we express
the second term as follows:
\begin{eqnarray}
&&\frac1tE_x(f(X_t)\ind_{\{T\le t,\,X_t\in I_{\delta,x}^c\cap(0,\infty)\}})=\nonumber\\
&&\quad\qquad\int_{I_{\delta,x}^c\cap(0,\infty)}f(y)P_x(t\ge
T\,|\,X_t=y)\frac{p_1(t^{-1/\alpha}(y-x))}{t^{1+1/\alpha}}\,dy\,.\label{326}
\end{eqnarray}
From (\ref{238}) and (\ref{239}), there is a constant $C_5>0$ such
that for $|x|$ sufficiently large,
\[p_1(x)\le C_5|x|^{-\alpha-1}\,,\]
hence there exist $C_6>0$ and $t_1$ (which may depend on $x$) such
that for all $y\in I_{\delta,x}^c$ and for all $0<t\le t_1$,
\[p_1(t^{-1/\alpha}(y-x))\le C_6 t^{1+1/\alpha}\,.\]
Therefore, from the Lebesgue theorem of dominated convergence and
Lemma \ref{lem1}, the expression in (\ref{326}) tends to 0 as $t$
goes to 0.

Now recalling our first equality, we have for any $\delta\in(0,x)$,
\begin{equation}\label{1492}
\lim_{t\rightarrow0}\frac1tE_x(f(X_t)\ind_{\{T\le
t,\,X_t\in(0,\infty)\}})=
\lim_{t\rightarrow0}\frac1tE_x(f(X_t)\ind_{\{T\le t,\,X_t\in
I_{\delta,x}\cap(0,\infty)\}})\,.\end{equation} Set
$b_{\delta,x}^-=\inf\{f(y),\,y\in I_{\delta,x}\}$ and
$b_{\delta,x}^+=\sup\{f(y),\,y\in I_{\delta,x}\}$. From our
hypothesis on $f$, $b_{\delta,x}^-$ and $b_{\delta,x}^+$ are finite
and from the equality above, we have
\begin{eqnarray}&&b_{\delta,x}^-\lim_{t\rightarrow0}\frac1t[P_x(T\le
t,\,X_t\in(0,\infty))-P_x(T\le t, X_t\in
I_{\delta,x}^c\cap(0,\infty))]\le\nonumber\\&&\;\;\;\lim_{t\rightarrow0}
\frac1tE_x(f(X_t)\ind_{\{T\le t,\,X_t\in(0,\infty)\}}) \le
\lim_{t\rightarrow0}b_{\delta,x}^+\frac1tP_x(T\le
t,\,X_t\in(0,\infty)).\label{5925}
\end{eqnarray}
But applying again (\ref{1492}) with $f\equiv1$, we find
$\lim_{t\rightarrow0}\frac1tP_x(T\le t, X_t\in
I_{\delta,x}^c\cap(0,\infty))=0$.

Now write:
\[\frac1tP_x(T\le t,\,X_t\in(0,\infty))=\frac1tP_x(T\le
t)-\frac1tP_x(X_t\in(-\infty,0])\,,\] and apply (\ref{437}) together with the fact
that $\lim_{t\rightarrow0}\frac1tP_x(X_t\in(-\infty,0])=c_-/(\alpha
x^\alpha)$, see for instance \cite{Be}, Exercise I.1.

Finally note that $\delta$ is arbitrarily small in the inequalities
(\ref{5925}) and since $f$ is continuous, $b_{\delta,x}^-$ and
$b_{\delta,x}^+$ tend to $f(x)$ as $\delta$ goes to 0. This allows
us to conclude.
\end{proof}

\section{Killed or conditioned stable processes as pssMp}

In this section, we compute the characteristics of the underlying
L\'evy process in the Lamperti representation of a pssMp $(X,\p_x)$
when this process is either a stable L\'evy process which is killed
when it first hits the positive half-line (section \ref{kil}) or a
stable L\'evy process conditioned to stay positive (section
\ref{pos}) or a stable L\'evy process conditioned to hit 0
continuously (section \ref{zero}). If $(X,P_x)$ is a stable
subordinator, then it can be considered as its own version
conditioned to stay positive and in this case and the
characteristics of the underlying L\'evy process have been computed
by Lamperti \cite{La}, Section 6. Except in this situation, the
cases where $(X,P_x)$ or $(-X,P_x)$ is a subordinator have no
interest in this study, so they will be implicitly excluded in the
sequel. Also, as already mentioned in the introduction, since all
our study is well known when $(X,P_x)$ is the standard Brownian
motion, we will always suppose that $\alpha\neq2$.

\subsection{The killed process}\label{kil}

In this subsection, we suppose that $(X,\p_x)$, $x>0$ is a stable
L\'evy process with index $\alpha \in (0,2)$ which is killed when it
first leaves the positive half-line. To define this process more
formally, let $(X,P_x)$ be a stable L\'evy process starting at
$x>0$. {\it We keep the same notations as in Section $\ref{prelim}$ for
the characteristics of $(X,P_x)$.} Recall that $T = \inf\{t\geq 0 :
X_t \leq 0\}$, then the probability measure $\p_x$ is the law under
$P_x$ of the process
\begin{equation}\label{defkil}
X_t\ind_{\{t<T\}}\,,\;\;\;t\ge0\,.\end{equation} (Note that rather
than the {\it killed} process, we could also call $(X,\p_x)$ the
initial L\'evy process $(X,P_x)$ {\it absorbed} at level 0). It is
not difficult to see, that the process $(X,\p_x)$ is a positive
self-similar Markov process with index $\alpha$ such that $S <
\infty$, $\p_x$-a.s. Furthermore, if $(X,P_x)$ has no negative
jumps, then $(X,\p_x)$ ends continuously at $0$, so it belongs to
the class $\mathcal{C}_2$. If $(X,P_x)$ has negative jumps, then it
is known that it crosses the level 0 for the first time by a jump,
so $(X,\p_x)$ ends by a jump at $0$ and belongs to the class
$\mathcal{C}_3$. We will compute the infinitesimal generator of
$(X,\p_x)$ and deduce from its expression the law of the underlying
L\'evy process $\xi$ associated to $(X,\p_x)$ in the Lamperti
representation.

Specializing  the expression given in the introduction for stable L\'evy
processes, we obtain the infinitesimal generator $\mathcal{A}$ with domain
$\mathfrak{D}_{\mathcal{A}}$ of the process $(X,P_x)$:
\begin{equation}\label{5469}
\mathcal{A}f(x) = af^{\prime}(x)
+ \int_{\R}(f(x+y)-f(x)-yf^{\prime}(x)\ind_{\{\vert y
\vert<1\}})\nu(y)dy\end{equation}
 for  $f\in \mathfrak{D}_{\mathcal{A}}$, where we recall from the beginning
 of  section
\ref{prelim} that $\nu(y) = c_+ y^{-\alpha-1}\textbf{1}_{\{y>0\}}
        + c_- \vert y \vert ^{-\alpha-1}\textbf{1}_{\{y<0\}}$
is the density of the the L\'evy measure and  that  $c_-\ge0$,
$c_+\ge0$, $a=\frac{c_+-c_-}{1-\alpha}$, if $\alpha\neq1$ and $a=0$,
$c_+=c_-$, if $\alpha=1$.

In the sequel, we will denote by $\mathcal{K}$ the infinitesimal
generator of the killed process $(X,\p_x)$. Note that since the
state space of this process is $[0,\infty)$ and 0 is an absorbing
state, the domain of $\mathcal{K}$, that will be denoted by
$\mathfrak{D}(\mathcal{K})$, is included in the set
$\{f:[0,\infty)\rightarrow\R:f(0)=0\}$. From the expression of the
infinitesimal generator $\mathcal{A}$, we can deduce this of
$\mathcal{K}$ as shows the following result.

\begin{theorem}\label{main}
Let $(X,\p_x)$ be the pssMp which is defined in $(\ref{defkil})$ and  let
$\mathcal{K}$ be its generator. Let $f\in\mathfrak{D}_{\mathcal{K}}$
such that the function $\tilde{f}$ defined on $\R$ by
$$\tilde{f}(x) = \left\{\begin{array}{ll}f(x)&\mbox{if $x>0$}\\
0&\mbox{if $x\le0$}\end{array}\right.,$$ belongs to
$\mathfrak{D}_\mathcal{A}$, then
\begin{eqnarray*}
\mathcal{K}f(x)=\mathcal{A}\tilde{f}(x)
-\frac{f(x)}{x^{\alpha}}\left(k-\frac{c-}{\alpha}\right), \,\
x>0\,,\;\;\;\mathcal{K}f(0)=0\,,
\end{eqnarray*}
where the constant $k$ is defined in Lemma $\ref{lem2}$. The
generator $\mathcal{K}$ can also be written as
\begin{eqnarray*}
\mathcal {K}f(x) &=&
\int_{\R^+}\frac{1}{x^{\alpha}}(f(ux)-f(x)-xf^{\prime}(x)(u-1)
\ind_{\{\vert u-1 \vert<1\}})\nu(u-1)du\\
&&\qquad\qquad+ a x^{1-\alpha}f^{\prime}(x) - k x^{-\alpha}f(x)\,.
\end{eqnarray*}
\end{theorem}

\noindent {\it Remark}: We emphasize that the set of functions which
is used in the above statement to describe the generator
$\mathcal{K}$ contains at least all functions of the set
$\{f:[0,\infty)\rightarrow\R:f(0)=0\}$ such that
$\tilde{f}\in\mathcal{C}_b^2(\R)$.

\begin{proof}
Recall that $T = \inf\{t\geq 0 : X_t \leq 0\}$, $S = \inf\{t\geq 0 :
X_t = 0\}$ and let $f$ be a function which is as in the statement of
the theorem. Then note that
\begin{eqnarray*}
\e_x(f(X_t))&=&\e_x(f(X_t)\textbf{1}_{\{t<S\}}+f(0)\ind_{\{t\ge
S\}})= E_x(\tilde{f}(X_t)\textbf{1}_{\{t<T\}})\\&=&
E_x(\tilde{f}(X_t))-E_x(\tilde{f}(X_t) \textbf{1}_{\{T\le t\}})\,.
\end{eqnarray*}
So, for any $x>0$ the generator of the killed process $(X,\p_x)$
is given by:
\begin{eqnarray*}
\mathcal {K}f(x)&=&\lim_{t
\to\infty}\frac{1}{t}\,\e_x(f(X_t)-f(x))\\
&=&\lim_{t\to 0}\frac{1}{t}\,[E_x(\tilde{f}(X_t))-\tilde{f}(x)] -
\lim_{t\to 0} \frac{1}{t}\,E_x (\tilde{f}(X_t)\textbf{1}_{\{
T\le t\}})\\
&=&\mathcal{A}\tilde{f}(x)- \frac{f(x)}{x^\alpha}\left(k-\frac{c-}{\alpha}\right)\,.\\
\end{eqnarray*}
Where the last equality comes from Lemma \ref{lem2}, since
$\tilde{f}$ is continuous and bounded. The value of $\mathcal{K}f$
at $0$ is easily computed.

To prove the second assertion of the theorem, write
\[\mathcal {K}f(x) = af^\prime(x) -
\frac{f(x)}{x^\alpha}(k-\frac{c-}{\alpha})
+\int_{\R}(\tilde{f}(x+y)-f(x)-yf^{\prime}(x)\ind_{\{\vert y
\vert<1\}})\nu(y)\,dy\]  and let $I$ be the integral term above.
Then make the change of variable $y =x(u-1)$ to obtain,
\[ I =\frac{1}{x^\alpha} \int_{u \in \R}[\tilde{f}(xu)-f(x)-x(u-1)f^{\prime}(x)
\ind_{\{\vert x(u-1) \vert<1\}}]\nu(u-1)\,du\,. \] We rewrite $I$ in
the following form:
\begin{eqnarray*}
I =& \frac{1}{x^{\alpha}} \int_{(u
>0)}[\tilde{f}(xu)-f(x)-x(u-1)f^{\prime}(x)\ind_{\{\vert
u-1 \vert<1\}}]\nu (u-1)\,du \\
  +&\frac{1}{x^{\alpha}} \int_{(u > 0)}[x(u-1)f^{\prime}(x)(\ind_{\{\vert u-1 \vert<1\}}
  - \ind_{\{\vert x(u-1)\vert<1\}}]\nu(u-1)\,du \\
  +&\frac{1}{x^{\alpha}} \int_{(u < 0)} [\tilde{f}(xu)-f(x)-x(u-1)f^{\prime}(x)
  \ind_{\{\vert x(u-1) \vert<1\}}]\nu(u-1)\,du
\end{eqnarray*}
and we call each if these integrals $I_1,I_2,I_3$ respectively.
Integral $I_1$ stays as it is but $I_2$ and $I_3$ require additional
calculations:
\[I_3 = -\frac{f(x)}{x^\alpha}\int_{(u<0)}\nu(u-1)du -\frac{1}{x^\alpha}
\int_{(u < 0)}x(u-1)f^{\prime}(x)\ind_{\{\vert x(u-1)
\vert<1\}}\nu(u-1)\,du\,.\] Now suppose that $\alpha\neq1$ (the case
$\alpha=1$ being much simpler). We may verify (after fastidious
calculations) that the sum of $I_2$ and  the second term of $I_3$
gives
\[\frac{c_+ -c_-}{1-\alpha}
(1-x^{\alpha-1})\frac{f^{\prime}(x)}{x^{\alpha - 1}} =
 af^{\prime}(x)( x^{1-\alpha}-1)\,,\]
since $a = \frac{c_+ -c_-}{1-\alpha}$. We finally calculate the
first term of $I_3$:
$$-\frac{f(x)}{x^\alpha}\int_{(u<0)}\nu(u-1)\,du
= - \frac{f(x)}{x^\alpha}\frac{c_-}{\alpha}\,.$$ Then by adding again all the different
parts together, we find for the expression of ${\cal K}f(x)$:
\begin{eqnarray*}
\mathcal{K}f(x)
&=&af^\prime(x) - \frac{f(x)}{x^\alpha}(k-\frac{c-}{\alpha}) + I_1 +
af^{\prime}(x)(x^{1-\alpha}-1) - \frac{f(x)}{x^\alpha}\frac{c-}{\alpha} \\
&=& \frac{a}{x^{\alpha-1}}f^{\prime}(x)+I_1 - \frac{f(x)}{x^\alpha}k\,,
\end{eqnarray*}
which ends the proof.
\end{proof}

Let $\xi$ be the underlying L\'evy process in the Lamperti
representation of $(X,\p_x)$, as it is stated in (\ref{lamp}).
Recall that $\xi$ may have finite lifetime, so its characteristic
exponent $\Phi$ is defined by
\begin{equation}\label{3129}
\e[\exp(i\lambda\xi_t)\ind_{\{t<\zeta(\xi)\}}]=\exp[t\Phi(\lambda)]\,,\;\;\;\lambda\in\R\,.
\end{equation}
Using Lamperti's result which is recalled in Theorem \ref{thlamp} in
the introduction and Theorem \ref{main}, we may now give the
explicit form of the generator of $\xi$ in the special setting of
this subsection.
\begin{coro}\label{cor1}
Let $\xi$ be the L\'evy process in Lamperti representation
$(\ref{lamp})$ of the pssMp which is $(X,\p_x)$ defined in
$(\ref{defkil})$. The infinitesimal generator $\mathcal{L}$ of $\xi$
with domain $\mathfrak{D}_{\mathcal{L}}$ is given by
$$\mathcal{L} f (x) = a f^{\prime}(x) + \int_{\R} (f(x+y) - f(x) -
f^{\prime}(x)(e^y-1)\ind_{\{\vert e^y-1 \vert<1\}})\pi(y)dy - k
f(x)\,,$$ for any $f\in\mathfrak{D}_{\mathcal{L}}$ and $x\in\R$,
where $\pi(y) = e^y\nu(e^y -1)$, $y\in\R$. Equivalently, the
characteristic exponent of $\xi$ is given by
\[ \Phi(\lambda) = ia\lambda +\int_{\R} (e^{i\lambda y}-1-i\lambda
(e^y-1) \ind_{\{\vert e^y-1 \vert<1\}}) \pi(y)dy - k\,.\] The
process $(X,\p_x)$ belongs to the class $\mathcal{C}_3$ if $k>0$ and
it belongs to the class $\mathcal{C}_2$ if $k=0$. In the first case
the L\'evy process $\xi$ has finite lifetime with parameter $k$, in
the second case, it has infinite lifetime.
\end{coro}
\noindent It is rather unusual to see $l(y)=(e^y-1)\ind_{\{\vert
e^y-1 \vert<1\}}$ as the compensating function in the expression of
the infinitesimal generator or the characteristic exponent of a
L\'evy process. However, as noticed in the introduction, any
function $l$ such that $l(y)\sim y$, as $y\rightarrow0$ may be
chosen and the more classical function
$l(y)=y\ind_{\{\vert y \vert<1\}}$, would have the effect of
replacing the parameter $a$ by another one which expression is rather
complicated.\\

\noindent  Let us consider the unkilled version of $\xi$, i.e. the
L\'evy process $\tilde{\xi}$ with characteristic exponent \[
\tilde{\Phi}(\lambda) = ia\lambda +\int_{\R} (e^{i\lambda
y}-1-i\lambda (e^y-1) \ind_{\{\vert e^y-1 \vert<1\}}) \pi(y)dy\,.\]
A natural question is to know wether if the process $\tilde{\xi}$
oscillates, drifts to $-\infty$ or drifts to $+\infty$. Let us show
that the three situations may happen depending on the relative
values of $c_-$, $c_+$ and $\alpha$. From the expression of
$\tilde{\Phi}$, we see that $\tilde{\xi}$ is integrable and
\begin{eqnarray}&&\e(\tilde{\xi}_1)=-i\tilde{\Phi}'(0)
=a+c_+\left(\int_0^{\log 2}\frac{(1+y-e^y)e^y}{(e^y-1)^{\alpha+1}}\,dy+
\int_{\log
2}^\infty\frac{ye^y}{(e^y-1)^{\alpha+1}}\,dy\right)\nonumber\\&&\qquad\qquad+c_-
\int_{-\infty}^{0}\frac{(1+y-e^y)e^y}{(1-e^y)^{\alpha+1}}\,dy\,.\label{exp1}
\end{eqnarray}
(Here $\p$ can be any of the measures $\p_x$, $x>0$). On the one
hand, it is clear from the classification which is recalled in the
introduction that when $(X,P_x)$ has no negative jumps
(i.e.~$c_-=0$), then  the L\'evy process $\tilde{\xi}=\xi$ drifts
towards $-\infty$, so that
\begin{equation}\label{exp2}
\frac{c_+}{1-\alpha}+c_+\left(\int_0^{\log
2}\frac{(1+y-e^y)e^y}{(e^y-1)^{\alpha+1}}\,dy+ \int_{\log
2}^\infty\frac{ye^y}{(e^y-1)^{\alpha+1}}\,dy\right)<0\,,\end{equation}
for all $c_+>0$ and $\alpha\in(1,2)$. (Recall that in the spectrally
one side case, we have necessarily $\alpha\in(1,2)$.) On the other
hand, when $(X,P_x)$ has no positive jumps ($c_+=0$), then it is
easy to derive from (\ref{exp1}) that for any $c_->0$ fixed,
$\lim_{\alpha\downarrow1}\e(\tilde{\xi_1})=+\infty$ and
$\lim_{\alpha\uparrow2}\e(\tilde{\xi_1})=-\infty$. Since
$\alpha\mapsto\e(\tilde{\xi}_1)$ is continuous, there are values of
$\alpha\in(1,2)$ for which $\tilde{\xi}$ drifts to $-\infty$,
oscillates or drifts to $+\infty$. This argument and (\ref{exp2})
show that for all $c_->0$ and $c_+>0$, there are values of
$\alpha\in(1,2)$ for which $\tilde{\xi}$ drifts to $-\infty$.

\subsection{The process conditioned to stay positive}\label{pos}

We consider again a stable L\'evy process $(X,P_x)$ as it is defined as in
section \ref{prelim}. Formally, the process $(X,P_x)$ conditioned to stay
positive is an $h$-transform of the killed process defined in section \ref{kil},
i.e.
\begin{equation}\label{condpos}
\p_x^\uparrow(A)=h^{-1}(x)E_x(h(X_t)\ind_A\ind_{\{t<T\}})\,,\;\;\;
x>0,\,t\ge0,\,A\in{\cal F}_t\,,
\end{equation}
where $h(x)=x^{\alpha\rho}$. The function $h$ being positive and
harmonic for the killed process, formula (\ref{condpos}) defines the
law of a strong homogeneous Markov process. Moreover this process is
$(0,\infty)$-valued and it is clear that it inherits the scaling
property with index $\alpha$ from $(X,P_x)$. Hence
$(X,\p_x^\uparrow)$ yields an example of pssMp which belongs to the
class $\mathcal{C}_1$. The following more intuitive (but no less
rigorous) construction of the law $\p_x^\uparrow$ justifies that
$(X,\p_x^\uparrow)$ is called {\it the L\'evy process $(X,P_x)$
conditioned to stay positive}
\[\p_x^\uparrow(A)=\lim_{t\rightarrow+\infty}P_x(A\,|\,T>t)\,,\;\;\;
x>0,\,t\ge0,\,A\in{\cal F}_t\,.\]
We refer to \cite{Ch} for
a general account on L\'evy processes conditioned to stay positive. In particular it is proved
in \cite{Ch} that $(X,\p_x^\uparrow)$ drifts to $+\infty$ as $t\rightarrow+\infty$, i.e.
\begin{equation}\label{drift}
\p_x^\uparrow\big(\lim_{t\rightarrow+\infty}X_t=+\infty\big)=1\,.
\end{equation}
Let us also mention that this conditioning has a discrete time counterpart for random
walks. Let $\mu$ be a law which is in the domain of attraction of the stable
law $(X_1,P_0)$ and let $S^\uparrow$ be a random walk with law $\mu$ which is conditioned to stay
positive. Then the process $(X,\p_x^\uparrow)$ may be obtained as the limit in law of the
process $(n^{1/\alpha}S_{[nt]}^\uparrow,\,t\ge0)$, as $n$ tends to $\infty$. This invariance principle
has recently been  proved in \cite{CC}.\\

Since $(X,\p_x^\uparrow)$ is an $h$-process of the killed process $(X,\p_x)$ defined at the
previous subsection, its infinitesimal generator, that we denote by $\mathcal{K}^\uparrow$,
may be derived from $\mathcal{K}$ as follows
\begin{equation}\label{hgen}
\mathcal{K}^{\uparrow} f(x) = \frac{1}{h(x)} \mathcal{K}
(hf)(x)\,,\;\;\;x>0\,,\;\; f\in
\mathfrak{D}_{\mathcal{K}^{\uparrow}}\,.\end{equation} From
(\ref{hgen}) and Theorem \ref{main}, we obtain for $x>0$ and $f\in
\mathfrak{D}_{\mathcal{K}^{\uparrow}}$:
\begin{eqnarray*}
x^{\alpha}\mathcal {K}^\uparrow f(x)&=&\frac1{x^{\alpha\rho}}
\int_{\R^+}[(hf)(ux)-(hf)(x)-x(hf)^{\prime}(x)(u-1)
\ind_{\{\vert u-1 \vert<1\}}]\nu(u-1)du\\
&&\qquad\qquad+ a x(hf)^{\prime}(x) - k
(hf)(x)\\
&=&\int_{\R^+}[u^{\alpha\rho}f(ux)-f(x)-(\alpha\rho
f(x)+xf^{\prime}(x))(u-1)
\ind_{\{\vert u-1 \vert<1\}}]\nu(u-1)du\\
&&\qquad\qquad+ axf^{\prime}(x)+(a\alpha\rho-k) f(x)\,.
\end{eqnarray*}
Let us denote by $J$ the integral term in the above expression and
define $\nu^\uparrow(u)=u^{\alpha\rho}\nu(u-1)$. Then
\begin{eqnarray*}
J&=&\int_{\R^+}[f(ux)-u^{-\alpha\rho}f(x)-(\alpha\rho
f(x)+xf^{\prime}(x))u^{-\alpha\rho}(u-1)\ind_{\{\vert u-1 \vert<1\}}]\nu^\uparrow(u)du\\
&=&\int_{\R^+}[f(ux)-f(x)-xf^{\prime}(x)(u-1)\ind_{\{\vert u-1
\vert<1\}}]\nu^\uparrow(u)du+\\
&&\int_{\R^+}[u^{\alpha\rho}-1-\alpha\rho (u-1)\ind_{\{\vert u-1
\vert<1\}}]\nu(u-1)du\,f(x)+\\
&&\int_{\R^+}(u^{\alpha\rho}-1)(u-1)\ind_{\{\vert u-1
\vert<1\}}\nu(u-1)du\,xf'(x)\,.
\end{eqnarray*}
The infinitesimal generator of the process $(X,\p_x^\uparrow)$  is
then
\begin{eqnarray}
\mathcal{K}^{\uparrow} f(x)&=&
\frac{1}{x^{\alpha}}\int_{\R^+}[f(ux)-f(x)-xf^{\prime}(x)(u-1)\ind_{\{\vert
u-1\vert<1\}}]\nu^\uparrow(u)du\nonumber\\ &&\qquad\qquad+(a+a_1)
x^{1-\alpha}f^{\prime}(x)+(a\alpha\rho+a_2-k)x^{-\alpha}f(x)\,,\label{genpos}
\end{eqnarray}
where
\begin{eqnarray}
\nu^\uparrow(u) &=& u^{\alpha\rho} \nu(u-1)\nonumber\\
a_1&=&c_+\int_0^1\frac{(1+x)^{\alpha\rho}-1}{x^\alpha}\,dx+c_-
\int_0^1\frac{(1-x)^{\alpha\rho}-1}{x^\alpha}\,dx\label{integre}\\
a_2&=&c_+\left(\int_0^1\frac{(1+x)^{\alpha\rho}-1-\alpha\rho
x}{x^{\alpha+1}}\,dx+
\int_1^\infty\frac{(1+x)^{\alpha\rho}-1}{x^{\alpha+1}}\,dx\right)\nonumber\\
&&+c_-\int_0^1\frac{(1-x)^{\alpha\rho}-1+\alpha\rho x}{x^{\alpha+1}}\,dx\label{integ}
\end{eqnarray}
and $a$, $k$ and $\rho$ are given in  section \ref{prelim}. Note
that since $(X,\p_x^\uparrow)$ belongs to ${\cal C}_1$, the killing
rate $a\alpha\rho+a_2-k$ in the expression (\ref{genpos}) of its
generator must be zero, which gives the following expression for the
constant $k$:
\begin{equation}\label{constk}
k=a\alpha\rho+a_2\,.
\end{equation}
From (\ref{1257}), the value of $k$ is explicit in terms of the
constants $c_+$, $c_-$ and $\alpha$, so we should be able to compute
the integrals in the expressions (\ref{integre}) and (\ref{integ})
of $a_1$ and $a_2$. However, the calculation of these integrals
relies to special functions and its seems to be possible to check
that (\ref{1257}) and (\ref{constk}) coincide only in the trivial
cases $\alpha\rho=1$ and $\alpha=1$.

As for $(X,\p_x)$ in the previous subsection, we may now apply
Theorem \ref{thlamp} together with (\ref{genpos}) to compute the
characteristics of the underlying L\'evy process in the Lamperti
representation of $(X,\p_x^\uparrow)$.
\begin{coro}\label{cor2}
Let $\xi^\uparrow$ be the L\'evy process in the Lamperti
representation $(\ref{lamp})$ of the pssMp $(X,\p_x^\uparrow)$ which
is  defined in $(\ref{condpos})$. The infinitesimal generator
$\mathcal{L}^\uparrow$ of $\xi^\uparrow$ with domain
$\mathfrak{D}_{\mathcal{L}^\uparrow}$ is given by
$$\mathcal{L}^\uparrow f (x) = a^\uparrow f^{\prime}(x) + \int_{\R} (f(x+y) - f(x) -
f^{\prime}(x)(e^y-1)\ind_{\{\vert e^y-1
\vert<1\}})\pi^\uparrow(y)dy\,,$$ for any
$f\in\mathfrak{D}_{\mathcal{L}^\uparrow}$ and $x>0$, where
 $\pi^\uparrow(y) = e^{(\alpha\rho+1)y}\nu(e^y
-1)$, $y\in\R$ and $a^\uparrow=a+a_1$, the constant $a_1$ being
defined in $(\ref{integre})$. Equivalently, the characteristic
exponent of $\xi^\uparrow$ is given by
\[ \Phi^\uparrow(\lambda) = ia^\uparrow\lambda +\int_{\R} (e^{i\lambda y}-1-i\lambda
(e^y-1) \ind_{\{\vert e^y-1 \vert<1\}}) \pi^\uparrow(y)dy\,.\]
\end{coro}
\noindent It follows from (\ref{drift}) that  the
underlying L\'evy process $\xi^\uparrow$ drifts to $+\infty$. This process being integrable,
it means in particular that $0<\e(\xi^\uparrow_1)=-i\Phi^\uparrow(0)<\infty$.

\subsection{The process conditioned to hit 0 continuously}\label{zero}

Let $S$ be an integer valued random walk which law
is in the domain of attraction of the stable law $(X_1,P_0)$. For
 $y\in\mathbb{Z}\setminus\{0\}$, define the law of the chain
$S^\searrow_y$ as this of the random walk $S_y$ starting from $y$ and
conditioned to hit 0 as follows~:
\begin{eqnarray*}&&(S^\searrow_y(n),\,0\le n\le\tau_{(-\infty,0]}^\searrow)
\ed[(S_y(n),\,0\le n\le\tau_{(-\infty,0]})\,|\,S_y(\tau_{(-\infty,0]})=0]\\
&&\qquad(S^\searrow_y(n),\,n\ge\tau_{(-\infty,0]}^\searrow)\equiv 0\,,
\end{eqnarray*}
where $\tau^\searrow_{(-\infty,0]}:=\inf\{n:S^\searrow_n\le0\}$ and
$\tau_{(-\infty,0]}:=\inf\{n:S_n\le0\}$. It has recently been proved
in \cite{CC}, that the rescaled linear interpolation of
$S_y^\searrow$, i.e.
 $$(n^{-1/\alpha}S^\searrow_{[n^{1/\alpha}x]}([nt]),\,t\ge0)\,,$$ converges in law
on the Skorohod's space as $n$ tends to $\infty$ towards a Markov
process which we will call here the L\'evy process $(X,P_x)$
conditioned to hit 0 continuously. Again, this process may be
defined more formally as an $h$-process of the killed process
$(X,\p_x)$ introduced in section \ref{kil}. In this case, the
positive harmonic function related to $(X,\p_x)$ is
$g(x)=x^{\alpha\rho-1}$, and for $x>0$, the law $\p_x^\searrow$ of
the conditioned process is defined by
\begin{eqnarray}\label{defzero}
&&\p_x^\searrow(A,t<S)={g(x)^{-1}}E_x(g(X_t)\ind_A\ind_{\{t<T\}})\label{hitcont}\\
&&\p_x^\searrow(X_t=0,\,\mbox{for all $t\ge S$})=1\,,\nonumber
\end{eqnarray}
for all $x>0$, $t\ge0$, and $A\in{\cal F}_t$.
It is proved in \cite{Ch} that the process $(X,\p_x^\searrow)$ reaches 0 continuously
(it may happen by an accumulation of negative jumps if $(X,P_x)$ has negative jumps), that is
\[\p_x^\searrow(X_{S-}=0)=1\,,\]
hence $(X,\p_x^\searrow)$ is a pssMp which belongs to the
class ${\cal C}_2$. The infinitesimal
generator of $(X,\p_x^\searrow)$ is
\begin{equation}\label{ggen}
\mathcal{K}^\searrow f(x) = \frac1{g(x)}\mathcal{K}
(gf)(x)\,,\;\;\;x>0\,,\;\; f\in
\mathfrak{D}_{\mathcal{K}^{\searrow}}\,.\end{equation} Trivially,
when there are no negative jumps (i.e.~$\alpha\rho=1$), $g\equiv1$
and the processes $(X,\p_x)$ and $(X,\p_x^\searrow)$ are the same.
The same calculations as in the subsection \ref{pos}, replacing
$\alpha\rho$ by $\alpha\rho-1$ lead to
\begin{eqnarray*}
\mathcal{K}^{\searrow} f(x)&=&
\frac{1}{x^{\alpha}}\int_{\R^+}[f(ux)-f(x)-xf^{\prime}(x)(u-1)\ind_{\{\vert
u-1\vert<1\}}]\nu^{\searrow}(u)du\\
&&\qquad+(a+a_3)
x^{1-\alpha}f^{\prime}(x)+(a(\alpha\rho-1)+a_4-k)x^{-\alpha}f(x)\,,
\end{eqnarray*}
where
\begin{eqnarray}
\nu^\searrow(u) &=& u^{\alpha\rho-1} \nu(u-1)\nonumber\\
a_3&=&c_+\int_0^1\frac{(1+x)^{\alpha\rho-1}-1}{x^\alpha}\,dx+c_-
\int_0^1\frac{(1-x)^{\alpha\rho-1}-1}{x^\alpha}\,dx\label{integre2}\\
a_4&=&c_+\left(\int_0^1\frac{(1+x)^{\alpha\rho-1}-1-(\alpha\rho-1)
x}{x^{\alpha+1}}\,dx+\int_1^\infty\frac{(1+x)^{\alpha\rho-1}-1}
{x^{\alpha+1}}\,dx\right)\nonumber\\
&&+c_-\int_0^1\frac{(1-x)^{\alpha\rho-1}-1+(\alpha\rho-1)
x}{x^{\alpha+1}}\,dx\,.\label{inte}
\end{eqnarray}
Here again, since $(X,\p_x^\searrow)$ belongs to ${\cal C}_2$, the
killing rate $a\alpha\rho+a_4-k$ of its generator must be zero,
which gives the following expression for the constant $k$:
\begin{equation}\label{constk2}
k=a(\alpha\rho-1)+a_4\,.
\end{equation}
Comparing (\ref{constk}) with (\ref{constk2}) we should be able to
check that
\begin{eqnarray*}
&&a_4-a_2=a=\frac{c_+-c_-}{1-\alpha}=\\&&
c_+\left(\int_0^1\frac{1-(1+x)^{\alpha\rho-1}}{x^\alpha}\,dx-
\int_1^\infty\frac{(1+x)^{\alpha\rho-1}}{x^\alpha}\,dx\right)+
c_-\int_0^1\frac{(1-x)^{\alpha\rho-1}-1}{x^\alpha}\,dx\,,\end{eqnarray*}
however, this seems to be possible to realise only in the trivial
cases $\alpha\rho=1$ and $\alpha=1$, $\rho=1/2$.

 As in the previous sections, we may now compute the
characteristics of the underlying L\'evy process in the Lamperti
representation of $(X,\p_x^\searrow)$.
\begin{coro}\label{cor3}
Let $\xi^\searrow$ be the L\'evy process in the Lamperti
representation $(\ref{lamp})$ of the pssMp $(X,\p_x^\searrow)$ which
is  defined in $(\ref{defzero})$. The infinitesimal generator
$\mathcal{L}^\searrow$ of $\xi^\searrow$ with domain
$\mathfrak{D}_{\mathcal{L}^\searrow}$ is given by
$$\mathcal{L}^\searrow f (x) = a^\searrow f^{\prime}(x) + \int_{\R} (f(x+y) - f(x) -
f^{\prime}(x)(e^y-1)\ind_{\{\vert e^y-1
\vert<1\}})\pi^\searrow(y)dy\,,$$ for any
$f\in\mathfrak{D}_{\mathcal{L}^\searrow}$ and $x\in\R$, where
$\pi^\searrow(y) = e^{\alpha\rho y}\nu(e^y -1)$, $y\in\R$ and
$a^\searrow=a+a_3$, the constant $a_3$ being defined in
$(\ref{integre2})$. Equivalently, the characteristic exponent of
$\xi^\searrow$ is given by
\[ \Phi^\searrow(\lambda) = ia^\searrow\lambda +\int_{\R} (e^{i\lambda y}-1-i\lambda
(e^y-1) \ind_{\{\vert e^y-1 \vert<1\}}) \pi^\searrow(y)dy\,.\]
\end{coro}
\noindent As we already noticed, the process
$(X,\p_x^\searrow)$ belongs to the class $\mathcal{C}_2$, therefore the
underlying L\'evy process $\xi^\searrow$ drifts to $-\infty$,
 in particular, since this process is also integrable,
  we have $-\infty<\e(\xi^\searrow_1)=-i\Phi^\searrow(0)<0$.

\section{The minimum of $\xi$ up to an independent
exponential time.}

With the same notations for $(X,\p_x)$, $(X,P_x)$, and $\xi$ as in
section \ref{kil}, we suppose here that $(X,P_x)$ has negative
jumps, that is $\alpha\rho<1$ (which is also equivalent to $c_->0$).
Recall that the characteristics of $\xi$ have been computed in
Corollary \ref{cor1}. The first result of this section consists in
computing an explicit form of the law of the overall minimum of
$\xi$. Since $\xi$ has finite lifetime, it has the same law as a
L\'evy process, say $\tilde{\xi}$, with infinite lifetime and killed
at an independent exponential time with parameter $k$. Then let us
show how Lamperti representation together with classical results on
undershoots of subordinators allow us to compute the law of the
minimum of $\tilde{\xi}$ up to an independent exponential time with
parameter $k$. The latter is known as the spacial Wiener-Hopf factor
of the L\'evy process $\tilde{\xi}$, see \cite{Do}.

Set $\underline{X}=\inf_{s\le S}X_s$ and $\underline{\xi}=\inf_{s\le
{\rm\bf \zeta}}\xi_s$, where we recall from the introduction that
$S=\inf\{t:X_t=0\}$ and $\zeta:=\zeta(\xi)$ is the lifetime of
$\xi$. Then on the one hand, from the Lamperti representation
(\ref{lamp}), under $\p_x$, the processes $\underline{X}$ and
$\underline{\xi}$ are related as follows:
\begin{equation}\label{eq1}
\underline{X}=x\exp\underline{\xi}\,,\;\;\;\p_x-\mbox{a.s.}
\end{equation}
On the other hand, let $H$ be the downward ladder height process
associated to $(X,P_0)$, that is $H_t=-X_{\eta_t}$, where $\eta$ is
the right continuous inverse of the local time at 0 of the process,
$(X,P_0)$ reflected at its minimum, i.e. $(X-\underline{X},P_0)$. We
refer to \cite{Be}, Chap. VI, for a definition of ladder height
processes. It is easy to see the following identity:
\begin{equation}\label{eq3}
\underline{X}=x-H_{\nu(x)-}\,,\;\;\;\p_x-\mbox{a.s.}
\end{equation}
where $\nu(x)=\inf\{t:S_t>x\}$. In other words, $\underline{X}$
corresponds to the so-called undershoot of the subordinator $H$ at
level $x$. Since $H$ is a stable subordinator with index
$\alpha\rho$, the law of $\underline{X}$, and hence this of
$\underline{\xi}$, can be computed explicitly as shown in the next
proposition. In the sequel, $\p$ will be a reference probability
measure under which $\xi$ and $H$ have the laws described above.
\begin{prop}\label{prop1}
Recall that $\rho=P_0(X_1\le0)$ and let $\xi$ be the L\'evy process
which law is described in Corollary $\ref{cor1}$. For any
$\lambda>0$,
\begin{equation}
\e(e^{\lambda
\underline{\xi}})=\frac{\Gamma(\lambda+1-\alpha\rho)}{\Gamma(\lambda+1)
\Gamma(1-\alpha\rho)}\,.
\end{equation}
In other words,  $\exp\underline{\xi}$ is a Beta variable with
parameters $\alpha\rho$ and $1-\alpha\rho$, i.e.
$\exp\underline{\xi}$ has density function:
$\p(\exp\underline{\xi}\in
dt)=\beta(\alpha\rho,1-\alpha\rho)^{-1}t^{\alpha\rho-1}
(1-t)^{-\alpha\rho}\ind_{\{t\in[0,1]\}}\,dt$.
\end{prop}
\begin{proof} Recall that the L\'evy measure $\theta(dy)$ of $H$ and its
potential measure $U(dy)$ are given by:
$$\theta(dy)=c_1y^{-(\alpha\rho+1)}\textbf{1}_{\{y>0\}}\,dy\;\;\;\mbox{and}\;\;\;\;
\int_0^\infty e^{-\lambda y}\,U(dy)=c_2\lambda^{-\alpha\rho}\,,$$
where $c_1$ and $c_2$ are positive constants. Then from Proposition
2 of \cite{Be}, Chap. III,
\[\p(H_{\nu(x)-}\in dy)=\ind_{\{y\in[0,x]\}}\int_x^\infty
U(dy)\theta(dz-y)\,,\] from which we obtain for all $\lambda\ge0$
and $\mu\ge0$,
\begin{eqnarray*}
\int_0^\infty e^{-\mu
x}E(e^{-\lambda(x-H_{\nu(x)-})})\,dx&=&\int_0^\infty
e^{-(\lambda+\mu)x}\int_0^x e^{\lambda
y}\int_x^\infty U(dy)\,\theta(dz-y)\,dx\\
&=&\int_0^\infty
e^{-(\lambda+\mu)x}\int_0^x\frac{c_1}{\alpha\rho}e^{\lambda
y}(x-y)^{-\alpha\rho}U(dy)\,dx\\
&=&\frac{c_1c_2(\lambda+\mu)^{\alpha\rho-1}}{\alpha\rho\mu^{\alpha\rho}}
\Gamma(1-\alpha\rho)=\frac{(\lambda+\mu)^{\alpha\rho-1}}{\mu^{\alpha\rho}}\,.
\end{eqnarray*}
It means that if $\varsigma$ is exponentially distributed with
parameter $\mu$ and independent of $H$, then
$\varsigma-H_{\nu(\varsigma)-}$ is gamma distributed with parameters
$\mu$ and $1-\alpha\rho$, i.e.
\[\e\left(e^{-\lambda(\varsigma-H_{\nu(\varsigma)-})}\right)=
\left(\frac{\mu}{\lambda+\mu}\right)^{1-\alpha\rho}\,.\] Recall that
the moment of order $\lambda>0$ of the Gamma law with parameters
$\mu$ and $1-\alpha\rho$ is
$\Gamma(\lambda+1-\alpha\rho)/(\mu^\lambda\Gamma(1-\alpha\rho))$,
then thanks to (\ref{eq1}) and (\ref{eq3}), one has
\[\e(e^{\lambda\underline{\xi}})=\frac{\e(\gamma^{\lambda})}
{\e(\varsigma^{\lambda})}=\frac{\Gamma(\lambda+1-\alpha\rho)}
{\Gamma(\lambda+1)\Gamma(1-\alpha\rho)}\,,\] which is the moment of
order $\lambda$ of a Beta variable with parameters $\alpha\rho$ and
$1-\alpha\rho$.
\end{proof}
\noindent In view of the result of Proposition \ref{prop1}, one is
tempted to compute the law of the overall minimum $\inf_{t\le{\rm\bf
e}(\mu)}\tilde{\xi}_t$ of the unkilled process $\tilde{\xi}$ before
an independent exponential time of {\it any} parameter $\mu>0$.
However although the pssMp which is obtained from
$(\tilde{\xi}_t,\,t\le{\rm\bf e(\mu)})$ through Lamperti
representation is absolutely continuous with respect to $(X,\p_x)$,
its law is not sufficiently explicit to apply the same arguments a
in Proposition \ref{prop1}.\\

We can still apply the same arguments as above to determine to law
of the overall minimum of the L\'evy process $\xi^\uparrow$ which is
defined in section \ref{pos}. Indeed, as we observed in this
section, $\xi^\uparrow$ drifts to $+\infty$, as well as the pssMp
$(X,\p_x^\uparrow)$, and from Lamperti representation the relation
\begin{equation}\label{eq12}
\underline{X}=x\exp\underline{\xi}^\uparrow\,,\;\;\;\p_x^\uparrow-\mbox{a.s.}
\end{equation}
holds. Moreover, the law of $(\underline{X},\p_x^\uparrow)$ is
explicit and may be found in \cite{Ch}, Theorem 5: for all $x>0$,
\[\p_x^\uparrow(\underline{X}\le
y)=\frac{x^{\alpha\rho}-(x-y)^{\alpha\rho}\ind_{\{y\le
x\}}}{x^{\alpha\rho}}\,.\] This allows us to state:
\begin{prop}\label{ruin}
Let $\xi^\uparrow$ be the L\'evy process which law is described in
Corollary $\ref{cor2}$. The law of the overall minimum
$\underline{\xi}^\uparrow$ of $\xi^\uparrow$ is given by:
\begin{equation}
\p(-\underline{\xi}^\uparrow\le
z)=(1-e^{-z})^{\alpha\rho}\ind_{\{z\ge0\}}\,.
\end{equation}
\end{prop}
\noindent This computation is closely related to risk theory and in
particular proposition \ref{ruin} provides an explicit form of the
ruin probability at level $z\ge0$, i.e. \[\p(\exists\, t\ge0,\,
z+\xi^\uparrow_t\le0)=\p(\underline{\xi}^\uparrow\le-z)=1-(1-e^{-z})\]
for this class of L\'evy processes, see the recent paper by Lewis
and Mordecki \cite{LM}.

\vspace*{.3in}

\noindent {\bf Acknowledgments}\\

This work has been partly done during the visit of the second author
at the University of Mexico UNAM. He his very grateful to this
University for its support. We thank Pat Fitzsimmons for fruitful
discussions and valuable comments.


\end{document}